\def\fnote#1{\footnote}

\newtheorem{lem}{Lemma}
\newtheorem{teo}{Theorem}
\newtheorem{pro}{Proposition}
\newtheorem{defin}{Definition}
\newcommand{\be}{\begin{equation}}

\newcommand{\ee}{\end{equation}}
\newcommand{\beqa}{\begin{eqnarray}}
\newcommand{\eeqa}{\end{eqnarray}}
\newcommand{\ba}{\begin{array}}
\newcommand{\bal}{\begin{array}{l}}
\newcommand{\ea}{\end{array}}
\newcommand{\bt}{\begin{teo}}
\newcommand{\et}{\end{teo}}

\newcommand{\vv}{{\bf V}}
\newcommand{\la}{\langle}
\newcommand{\ra}{\rangle}

\documentclass[jsco]{academic}                                 

\input psfig.sty

\newcommand{\putfig}[5]{
\begin{figure}[hptb]
        \begin{center}
            \leavevmode
            \psfig{file=#1,height=#2,width=#3}
        \end{center}
        \vspace{-2ex}
        \caption{#4}
        \label{#5}
\end{figure}
}

\date{December 6, 2001}
\begin{document}
\shortauthor{Jarrah, Laubenbacher, Romanovski}
\shorttitle{The center variety}

\title{ The Center Variety of Polynomial Differential Systems}

\author[1]{Abdul Salam Jarrah}
\author[1]{Reinhard Laubenbacher}
\author[2]{Valery Romanovski}

\address[1]{Department of Mathematical Sciences, New Mexico State University,
Las Cruces, NM 88003, USA}
\address[2]{Center for Applied Mathematics and Theoretical Physics,
University of Maribor, Krekova 2, SI-2000 Maribor, Slovenia}

\maketitle

{\abstract{
We investigate the symmetry component of the center variety
of polynomial differential systems, corresponding to systems with an axis
of symmetry in the real plane. We give a general algorithm to find this
irreducible subvariety and compute its dimension.
We show that our methods provide a simple way to compute the radical of the ideal
generated by the focus quantities and, therefore, to estimate the cyclicity
of a center in the case when the ideal is radical.  In particular,
we use our methods
to get a simple  proof of the famous Bautin theorem on the cyclicity of
the  quadratic system.
}}

\section{Introduction}
We start with a brief history and background. For more information and proofs
for the facts we state without proofs, see either
the indicated original source or \cite{ALS, RosR, RoSc, Dana1, Dana2, Z}.

In his M\'emoire \cite{P}, Poincar\'e initiated the study of
dynamical systems by studying real polynomial differential systems
of the form
\begin{eqnarray} \label{mg}
\frac{du}{dt} =  U(u,v),\\ \nonumber
\frac{dv}{dt} =  V(u,v),
\end{eqnarray}
where $U$ and $V$ are polynomials over the real numbers with
a singular point $(u_0,v_0)$,
taken to be $(0,0)$, without loss of generality.
There he also defined many now standard concepts, in particular,
the notion of a {\it center} of a system. The origin is a center if there exists
a neighborhood $\cal U$ of $(0,0)$ such that every point of $\cal U$
other than $(0,0)$ is nonsingular, and the integral curve passing
through that point is closed. Moreover, he proved the following theorem.

\begin{teo}
Assume that the linearization of system (\ref{mg}) at the origin has purely
imaginary eigenvalues.  Then, without loss of generality, it is of the form
\begin{eqnarray}\label{ruv}
\frac{du}{dt} &= & -u+U_1(u,v) =  U(u,v),\\ \nonumber
\frac{dv}{dt} &= & v+V_1(u,v) =  V(u,v),
\end{eqnarray}
where $U_1$ and $V_1$ are polynomials with nonlinear terms.
Then the origin is a center if and only if there exists
a formal power series $\Psi(u,v)= u^2+v^2+\dots $,
convergent in a neighborhood of the origin, such that
\[
  {{\partial   \Psi}\over   {\partial   u}}U+
{{\partial   \Psi}\over {\partial v}}V \equiv 0.
\]
\end{teo}
The function $\Psi(u,v)$ is a local first integral of system
(\ref{mg}).

Lyapunov \cite{Lyap} generalized and proved the above theorem for the case
when $P$ and $Q$ are real analytic functions.

There are many different ways to enclose the set of plane real system (RS),
into the set of two dimensional complex systems (CS).
The  most  convenient and commonly used way is the following.
Consider the real plane $(u,v)$ as the complex line with the variable $x$:
\be \label{cn}
x=u+iv.
\ee
Then System (\ref{mg}) is equivalent to the equation
\be \label{cxpiy}
i\frac{dx}{dt}=P(x,\bar x),
\ee
where $P=iU-V$. 
In many cases it is convenient to use  just equation (\ref{cxpiy}), 
however it is natural to add to this equation its complex conjugate,
$- i\dot {\bar x}=\bar P(x,\bar x)$ and  consider $y=\bar x$ as a
new variable and $ Q=\bar P$ as a new function.
As a result, we get the system of two complex differential equations
\be \label{ns}
i\dot x =P(x,y), \ \ - i\dot {y}=Q(x,y).
\ee
Thus there is one-to-one correspondence between
 systems  from (RS)  and the subset of  (CS) consisting of
 systems  of the form (\ref{ns}) where the second equation is the complex
conjugate to the first one.

For  polynomial systems of the form  (\ref{ruv}) the procedure above
yields the system
\begin{eqnarray}\label{mg2}
i\frac{dx}{dt} &= & x+P_1(x,y) =  P(x,y),\\ \nonumber
-i\frac{dy}{dt} &= &  y+Q_1(x,y) =  Q(x,y),
\end{eqnarray}
where $P_1$ and $Q_1$ are complex polynomials with nonlinear terms
and $Q_1 = \overline{P_1}$.
After the change of time $id\tau=dt$ we can write (\ref{mg2}) in the form
\begin{eqnarray}\label{mmg2}
 \frac{dx}{d\tau} &= & x+P_1(x,y) =  P(x,y),\\ \nonumber
-\frac{dy}{d\tau} &= & y+Q_1(x,y) =  Q(x,y).
\end{eqnarray}

In \cite{Dul},  Dulac considered  system (\ref{mmg2}) where 
$P_1$ and $Q_1$ are arbitrary complex polynomials with nonlinear terms.
Moreover, he gave the following definition for a center at the origin.

 \begin{defin}
System (\ref{mmg2}) has a center at the origin if 
there is  an analytic first integral of the form
\be \label{Int}
\Psi(x,y)= xy+
\sum_{s=3}^{\infty} \sum_{j = 0}^s v_{j,s-j} x^j y^{s-j},
\ee
\end{defin}
where the $v_{j,s-j}$ are functions in the coefficients of $P$ and $Q$.

\smallskip
If (\ref{mmg2}) is the  complexification of (\ref{ruv}) by means of (\ref{cn})
then this definition is in agreement with the definition given by Poincar\'e.

When $P$ and $Q$ are quadratic polynomials,
Dulac gave necessary and sufficient conditions on the coefficients of
$P$ and $Q$ such that system (\ref{mmg2}) has a center at the origin.
Moreover, he asked if one can find necessary and sufficient conditions on the
coefficients of $P$ and $Q$ (of any given degree) such that system (\ref{mmg2})
has a center at the origin. This is the  so-called  {\it center problem}.

In this paper, we present
a partial solution to the center problem for  polynomial systems,
using methods from computational algebra.

Any polynomial system of the form (\ref{mmg2}) can be written in the form
\begin{eqnarray} \label{gs}
\frac{dx}{d{\tau}} &= & x-\sum_{(p,q) \in S}
a_{pq}x^{p+1}{y}^{q}= P(x,y),\\ \nonumber
-\frac{dy}{d{\tau}} &= & y-\sum_{(p,q) \in S}
b_{qp}x^{q}{y}^{p+1}= -Q(x,y),
\end{eqnarray}
where $P(x,y),Q(x,y)\in\mathbb{C}[x,y]$ and
$$
S =  \{(p_i,q_i) \: | \: p_i+q_i\geq 1,\: i= 1,\ldots ,l\}
\subset \{\{-1\}\cup \mathbb{N}\}\times \mathbb{N}.
$$
Throughout this paper, $\mathbb{N}$ is the set of nonnegative integers.

We denote by $E(a,b)(= \mathbb{C}^{2l})$ 
the parameter space of (\ref{gs}),
and by $\mathbb{C}[a,b]$ the polynomial ring in the variables
$a_{pq},b_{qp}$.

As we have shown above, in the case when
\be \label{cond}
y= \overline x,\, b_{ij}= \overline{a}_{ji},\, i d\tau= dt
\ee
system (\ref{gs}) is equivalent to
the system
\be \label{rgs}
i\frac{dx}{dt}= x-\sum_{(p,q) \in S}
a_{pq}x^{p+1}{\bar x}^{q},
\ee
which has  a center or focus at the origin in the real plane
$\{(u,v)\mid x= u+iv\}$,
where the system can be also  written in the form
(\ref{ruv})
$$
\dot u =  -v+U_1(u,v), \ \dot v = u+V_1(u,v).
$$
In this case we denote the parameter space by $E(a)$.

\begin{defin}
We say that System (\ref{gs}) has a center on the set $W\subset E(a,b)$,
if for any point $(\tilde a,\tilde b)\in W$  the
corresponding system (\ref{gs}) has a center at  the origin.
\end{defin}

It is known \cite{Dul}
that one can always find a Lyapunov function $\Psi$
of the form (\ref{Int}) such that
\begin{eqnarray} \label{bar4}
\frac{\partial \Psi}{\partial x}P(x,y)+
\frac{\partial \Psi}{\partial y}
Q(x,y)=  g_{11}\cdot (xy)^{2}+{g}_{22}\cdot (xy)^{3}+{g}_{33}\cdot (xy)^{4}+\cdots ,
\end{eqnarray}
where the $g_{ii}$ are  polynomials of $ \mathbb{C} [a,b]$ called
{\it focus quantities}.
Thus, the maximal set $V\subset E(a,b)$, on which System (\ref{gs})
has a center, is the set where all polynomials $g_{ii}, i= 1,2,\dots$, vanish,
that is, $V$ is the variety of the ideal generated by the focus quantities
$g_{ii}$.

Denote by $\vv (I)$ the variety of the ideal $I$.
\begin{defin}
The set
$$
  V= \vv (\la g_{11}, g_{22}, \dots, g_{ii}, \dots \ra)
$$
is called the center variety of  system (\ref{gs}).
\end{defin}
So, for every point in $V$ the corresponding system has a center at the
origin in the sense that there is a first integral of the form (\ref{Int}).
However, if $(a,b)\in V$ and $a_{pq}= \bar b_{qp}$ for all $(p,q)\in S$,
then  such a point corresponds to a real system of the form (\ref{rgs}),
which then has a topological center at the origin in the plane $x= u+iv$.
(For a geometrical interpretation of the center of the complex system (\ref{gs})
see, e.g., \cite{Z}.)

Therefore, given a system of the form (\ref{gs}),
the problem of finding the center variety (center problem) of
 the system arises.

Among the components of the center variety  of polynomial system (\ref{gs})
there are at least two components which can be found without computing
any focus quantities: one component consists of Hamiltonian systems, and the
other one,
which we call the {\it symmetry component}, corresponds to systems which have
an axis of symmetry
in the real plane $x= u+iv$. The component of
Hamiltonian systems has a simple geometry
because it is equal to the   intersection of linear
subspaces of $E(a,b)$.

To find the symmetry components one can proceed as follows.
With   system (\ref{gs}) we associate the linear operator
\begin{eqnarray}\label{L}
 L(\nu) &=& {L^1(\nu) \choose L^2(\nu)} \\
&=& {p_1 \choose q_1} \nu_1+\cdots +{p_l \choose q_l} \nu_l
+{q_l \choose p_l}\nu_{l+1}+\cdots +{q_1 \choose p_1} \nu_{2l},
\nonumber 
\end{eqnarray}
where $ (p_m,q_m)\in S$.
Let $\cal M$ denote the set of
all solutions $\nu= (\nu_1,\nu_2,\dots,\nu_{2l})$
with non-negative components of the equation
\be \label{lk}
L(\nu)= \left(\begin{array}{cc} k\\k\end{array}\right),
\ee
for all $k\in \mathbb{N}$.
Obviously, $\cal M$ is an Abelian monoid.
For $k= 0$, finding the solutions of equation (\ref{lk}) is
just the standard integer programming problem, see, e.g., \cite[Sect. 1.4]{St2}.
Our algorithm is an adaptation of the algorithm to solve the standard
problem.
Let $\mathbb{C}[{\cal M}]$ denote the subalgebra of $\mathbb{C}[a,b]$
generated by all monomials of the form
$$
a_{p_1q_1}^{\nu_1}a_{p_2q_2}^{\nu_2}\cdots a_{p_lq_l}^{\nu_l}
b_{q_lp_l}^{\nu_{l+1}}b_{q_{l-1}p_{l-1}}^{\nu_{l+2}}\cdots b_{q_1p_1}^{\nu_{2l}},
$$
for all $\nu\in\cal M$.  In order to simplify notation we will abbreviate
such a monomial by $[\nu ]= [\nu_1,\ldots ,\nu_{2l}]$.  For $\nu\in\cal M$,
let
$$
\bar\nu= (\nu_{2l},\nu_{2l-1},\ldots ,\nu_1).
$$
Furthermore, let
$$
{\rm IM}[\nu]= [\nu]-[\bar\nu],\hspace{.1in}{\rm RE}[\nu]= [\nu]+[\bar\nu].
$$

It is shown in \cite{R93,RR01} that the focus quantities  of
system (\ref{gs}) belong to
$\mathbb{C}[\cal M]$ and have the  form
\begin{equation} \label{GGI}
g_{kk}= \sum_{L(\nu) =  (k,k)^t}g_{(\nu)}IM[\nu],
\end{equation}
with $g_{(\nu)}\in {\bf Q},\ k= 1,2,\dots$.
(Similar properties of the focus quantities were also obtained
in \cite{Isp,Liu2}.)

Consider the ideal
$$I_{sym}=\la [\nu]-[\bar \nu]\mid\nu\in\cal M\ra \subset \mathbb{C}
[\cal M],$$
called the {\it Sibirsky ideal} of (\ref{gs}).
The  next statement 
is an obvious  generalization of Sibirsky's symmetry
criteria \cite{Sib1} for a center for real systems
to systems of the form (\ref{gs}).
\begin{pro} \label {prop1}
The system (\ref{gs}) has a center on the set $\vv (I_{sym})$.
\end{pro}

\begin{defin} \label{dsc}  The set  $\vv  (I_{sym})$
is called the symmetry component of
the center
variety.
\end{defin}
For  points $(a,b)\in  \vv (I_{sym})$ for which $a_{pq}=\bar b_{qp}$
the corresponding systems (\ref{rgs})  have an axis of symmetry (see Section 3).
This is the reason for the name ``symmetry component.''
However we should mention that 
we do not know a proof that $\vv(I_{sym})$ is indeed 
a component of the center variety, that is, it is not a 
proper subvariety of an irreducible subvariety of the center variety.
However 
Proposition 1 implies that $\vv (I_{sym})$ is a  subset of the center variety
and  we shall  show below, that  $\vv  (I_{sym})$  is irreducible.
Moreover, for all cubic polynomial
systems    $\vv  (I_{sym})$ we have investigated, it is indeed a
component. Thus we conjecture that for any polynomial system of the form 
(\ref{gs}), $\vv (I_{sym})$ is a component 
of the center variety.

In the present paper we give a simple and effective
algorithm for  computing generators for the Sibirsky ideal,
using methods from computational algebra.  This allows us to find
a finite set of defining polynomials for the
symmetry component of the center variety.
Moreover, we prove the following results.
\begin{teo} \label{sip}
The Sibirsky ideal
$I_{sym}$ is prime
in  $ \mathbb{C}[a,b]$.
\end{teo}
\bt \label{siv}
The dimension of the symmetry component, $\vv (I_{sym})$,
of the center variety is equal to  $l$ if all coefficients on
the right-hand side of system (\ref{gs}) are resonant and
$l+1$ otherwise.
\et
Recall that the resonant coefficients are the ones
which cannot be canceled by transformation of the system
to the Poincar\'e-Lyapunov normal form, that is, the coefficients of the
form $a_{kk}, b_{kk}$.

It should be mentioned that a substantial amount of literature
is devoted to different particular subfamilies of polynomial systems, mainly
to systems of the second to fifth degree (the bibliography
on  quadratic systems alone by J.~Reyn \cite{Reyn} contains approximately
1500 references). However, the Symmetry Component
Algorithm in Section 2, and Theorems \ref{sip} and \ref{siv}
are among the very few statements known up to now about the whole
class of polynomial systems.

The center  problem is closely connected to the $cyclicity $
problem, which is sometimes called the local 16th Hilbert problem \cite{FY}.

\begin{defin}
Let   $n_{a,\epsilon}$
be the number of limit cycles of
system (\ref{rgs}) in an $\epsilon$-neighborhood of the origin. Then
 we say that the singular point $x=0$  of system (\ref{rgs})
with given coefficients  $a^*\in E(a)$ has cyclicity
 $k$ with respect to the space $ E(a)$, if there  exist $ \delta_0,
 \epsilon_0$ such that
for every $0<\epsilon<\epsilon_0 $ and $0<\delta<\delta_0$
$$\max_{a\in U_\delta (a^*)} n_{a,\epsilon} =k.$$
 \end{defin}
The cyclicity of the origin
of real quadratic systems (that is, system  (\ref{rgs}) with
quadratic non-linearities)
was first investigated by Bautin \cite{Bau},
and by Sibirsky \cite{Sibc} (see also  \cite{ZN}) for the system (\ref{rgs}) with
homogeneous cubic nonlinearities. They proved that the cyclicity of
the systems is less or equal two, respectively four.
(If we take into account perturbations by linear terms,
then the cyclicity
is 3 and 5, respectively).

 In this paper  we  show
that Theorem \ref{sip} along with other methods
from computational algebra provides an efficient tool to investigate  the cyclicity
of polynomial systems in those cases,  where the ideal of focus
quantities $\la g_{11}, g_{22},\dots \ra$ is radical.
Moreover, even when
the ideal of focus quantities is not radical, we use our algorithm along with
some standard methods to compute the cyclicity of some cubic systems, see
\cite{JLR}.

\section {An Algorithm for the Symmetry Component}

In this section we give an algorithm to find a finite set of
generators for the Sibirsky ideal, hence for the symmetry component
of the center variety.  It works for general systems
(\ref{gs}).
As a corollary, we obtain a Hilbert basis, that is, a finite minimal
generating set, for the monoid $\cal M$ described in the introduction.

Let
$$
{\cal{A}} =  \left[ \matrix{ c_{11} & c_{21} & \dots  & c_{l1} &
                         c_{l2} & \dots & c_{22} & c_{12} \cr
                        c_{12} & c_{22} & \dots & c_{l2} & c_{l1} & \dots &
                        c_{21} & c_{11} } \right]
$$
be a $(2 \times 2l)$-matrix with entries $c_{ij} \in {\mathbb{Z}}$.
For each $k \in \mathbb{N}$,
let 
\[
{\cal V}(k) = \left\{ \nu = (\nu_1, \dots , \nu_{2l}) \in {\mathbb{N}^{2l} \mid {\cal{A}}}
\cdot \nu^t = \left( \matrix{k \cr k }\right) 
\right\}.
\]

The proof of the following lemma is straightforward.

\begin{lem}
Let ${\cal{M_A}} = \cup_{k \in \mathbb{N}} {\cal V}(k)$. Then $\cal{M_A}$ is an
Abelian submonoid of $\mathbb{N}^{2l}$.
Also, if 
$$
\mu =(\mu_1,\ldots ,\mu_{2l})\in {\cal{M_A}},
$$
then $\bar{\mu}=(\mu_{2l},\ldots ,\mu_1) \in {\cal{M_A}}$.
\end{lem}

Let $d=2l$ and
$R = \mathbb{C}[x_1, \dots, x_d]$ be the polynomial ring in $d$
variables.
We consider the following binomial ideal in $R$:
$$
 I_{\cal A} = \langle {\bf x}^\mu - {\bf x}^{\bar{\mu}} \mid \mu \in
                  {\cal{M_A}}  \rangle .
$$

We will obtain a Hilbert basis for ${\cal{M_ A}}$ from a Gr\"obner basis
of the ideal $I_{\cal A}$.
If the matrix $\cal A$ arises from the coefficients of system (1)
as follows:
$$
\cal A=\left(\begin{array}{cccccccc}
p_1&p_2&\ldots &p_l&q_l&q_{l-1}&\ldots &q_1\\
q_1&q_2&\ldots &q_l&p_l&p_{l-1}&\ldots&p_1
\end{array}\right),
$$
then $I_{\cal A}$ is precisely the Sibirsky ideal of (\ref{gs}).

We will represent $I_{\cal A}$
as the kernel of a homomorphism of polynomial rings, so we can
use a standard algorithm to compute a finite generating set for
this ideal.
As a consequence, we obtain an
algorithm to compute a finite set of polynomials that define
the symmetry component of the center variety of (1).

Let $S = \mathbb{C}[t_1^{\pm},t_2^{\pm},y_1, \ldots, y_l]$.
Define a ring homomorphism $\phi : R \longrightarrow S$ by
\begin{eqnarray*}
        x_i &\longmapsto& y_it_1^{c_{i1}}t_2^{c_{i2}}\\
x_{l+i} &\longmapsto& y_{l-i+1}t_1^{c_{{l-i+1},2}}t_2^{c_{{l-i+1},1}}
\end{eqnarray*}
for $i = 1, \ldots, l$.

\begin{teo} \label{iai}
${\rm{ker}(\phi)}=I_{\cal A}$.
In particular, $I_{\cal A}$  is a prime ideal in $R$.
\end{teo}
\begin{proof}
The second statement follows immediately from the first,
since $S$ is a domain.

To prove that ${\rm{ker}}(\phi)=I_{\cal A}$, we first show that
${\rm ker}(\phi)$ is a binomial ideal.  We can factor the map $\phi$
as follows.  Let
$$
T=k[t_1,t_2,s_1,s_2,x_1,\ldots ,x_d,y_1,\ldots ,y_l],
$$
and
$$
\phi':\mathbb{C}[x_1,\ldots ,x_d]\longrightarrow T
$$
be defined like $\phi$, except that if $c_{i_j}<0$, then it appears
as an exponent of the variable $s_j$ instead of $t_j$.
Then $\phi$ is equal to the composition of $\phi'$ followed by
the projection
$$
T\longrightarrow T/\langle t_1s_1-1,t_2s_2-1\rangle .
$$
It is straightforward to see that ${\rm ker}(\phi)={\rm ker}(\phi')$. Let
\[
J=\langle x_i-y_it_1^{c_{i_1}}t_2^{c_{i_2}},
x_{l+i}-y_{l-i+1}t_1^{c_{l-i+1,2}}t_2^{c_{l-i+1,1}}\mid i=1,\ldots ,l\rangle \subset T.
\]
Then it follows immediately from \cite[Theorem 2.4.2]{AL}
that ${\rm ker}(\phi)={\rm ker}(\phi')=J\cap R$.
We obtain a generating set for $J\cap R$ by computing a reduced Gr\"obner
basis for $J$ using an elimination ordering with $x_i<y_j,t_k,s_r$
for all $i,j,k,r$, and then intersecting it with $R$.  Since $J$
is generated by binomials, any reduced Gr\"obner basis also consists of
binomials.  This shows that ${\rm ker}(\phi)$ is a binomial ideal.

Now, let ${\bf x}^{\alpha}-{\bf x}^{\beta}$ be a binomial in $R$.  We
may assume that the two monomials have no common factors, that is,
${\rm supp}(\alpha)\cap {\rm supp}(\beta)=\emptyset$.  Then
\begin{eqnarray}\label{E:ker} \nonumber
\phi({\bf x}^\alpha - {\bf x}^\beta) =
(y_1t_1^{c_{11}}t_2^{c_{12}})^{\alpha_1} \cdots
(y_lt_1^{c_{l1}}t_2^{c_{l2}})^{\alpha_l}(y_lt_1^{c_{l2}}t_2^{c_{l1}})^{\alpha_
{l+1}}
\cdots (y_1t_1^{c_{12}}t_2^{c_{11}})^{\alpha_d}                            \\
 - (y_1t_1^{c_{11}}t_2^{c_{12}})^{\beta_1} \cdots
(y_lt_1^{c_{l1}}t_2^{c_{l2}})^{\beta_l}(y_lt_1^{c_{l2}}t_2^{c_{l1}})^{\beta_{l
+1}}\nonumber
\cdots (y_1t_1^{c_{12}}t_2^{c_{11}})^{\beta_d}                              \\
 =y_1^{(\alpha_1+\alpha_{2l})} \cdots
y_l^{(\alpha_l+\alpha_{l+1})}
 t_1^{(c_{11}{\alpha_1}
+ \cdots +c_{l1}{\alpha_l}+c_{l2}{\alpha_{l+1}} + \cdots +
c_{12}{\alpha_d})}
t_2^{(c_{12}{\alpha_1} + \cdots + c_{l2}{\alpha_l}+c_{l1}{\alpha_{l+1}} +
\cdots +
c_{11}{\alpha_d})}\nonumber                                                          \\
-  y_1^{(\beta_1+\beta_{2l})} \cdots
y_l^{(\beta_l+\beta_{l+1})}
 t_1^{(c_{11}{\beta_1}
+ \cdots +c_{l1}{\beta_l}+c_{l2}{\beta_{l+1}} + \cdots + c_{12}{\beta_d)}}
t_2^{(c_{12}{\beta_1} + \cdots + c_{l2}{\beta_l}+c_{l1}{\beta_{l+1}} + \cdots
+
c_{11}{\beta_d})}\nonumber
\end{eqnarray}
Thus, $\phi({\bf x}^\alpha - {\bf x}^{\beta}) = 0$ if and only if
$$
\mbox{ for all } i \geq 1, \, \alpha_i+\alpha_{d-i+1} = \beta_i+\beta_{d-i+1},
$$
\begin{eqnarray*}
c_{11}{\alpha_1}+ \cdots +c_{l1}{\alpha_l}+c_{l2}{\alpha_{l+1}} 
&+& \cdots + c_{12}{\alpha_d} \\
 &=& c_{11}{\beta_1} + \cdots +c_{l1}{\beta_l}+c_{l2}{\beta_{l+1}} 
+ \cdots + c_{12}{\beta_d},
\end{eqnarray*}
and
\begin{eqnarray*}
c_{12}{\alpha_1} + \cdots + c_{l2}{\alpha_l}+c_{l1}{\alpha_{l+1}} 
&+& \cdots + c_{11}{\alpha_d} \\
&=& c_{12}{\beta_1} + \cdots + c_{l2}{\beta_l}+c_{l1}{\beta_{l+1}}
+ \cdots + c_{11}{\beta_d}.
\end{eqnarray*}
Since $\rm{supp}(\alpha)\cap\rm{supp}(\beta)=\emptyset$, we obtain
the following facts from the first condition in (\ref{E:ker}).
If, for some $i$, $\alpha_i=0=\beta_i$, then
$\alpha_{d-i+1}=0=\beta_{d-i+1}$, so that $\beta_i=\alpha_{d-i+1}$.
Now suppose that $\beta_i\neq 0$ for some $i$.  Then $\alpha_i=0$.
But then $\beta_{d-i+1}=0$, otherwise $\alpha_{d-i+1}=0$, which
cannot be since $\beta_i\geq 0$ for all $i$.  And this, in turn,
implies that $\beta_i=\alpha_{d-i+1}$.  Hence
$$
\beta=(\beta_1,\ldots ,\beta_d)=(\alpha_d,\ldots ,\alpha_1)=\overline{\alpha}.
$$
It follows from the last two equations in (\ref{E:ker}) that $\alpha\in\cal{M_A}$.
This completes the proof of the theorem.
\end{proof}

\begin{teo} \label{hb}
With notation as above, let $G$ be a reduced Gr\"obner basis of $I_{\cal A}$,
with respect to some term ordering on $R$.
Then
\begin{enumerate}
\item The set
$$
{\cal H} = \{ \mu, \overline{\mu} \mid {\bf x}^\mu - {\bf x}^{\overline{\mu}}
\in G\} \cup  \{ {\bf e}_i +{\bf e}_j \mid j = d-i+1; i =1,\ldots, l\},
$$
where ${\bf e}_i= (0,\dots,0,1,0,\dots,0)$ is the $i$th basis vector, 
is a Hilbert basis of $\cal{M_A}$.
Note that ${\bf e}_i+{\bf e}_{d-i+1}\in \cal{M_A}$ for all $i$.
\item $I_{\cal{A}} = \langle {\bf x}^\mu -{\bf x}^{\overline{\mu}}\mid\mu\in\cal H\rangle$.
\end{enumerate}
\end{teo}
\begin{proof}
The proof is similar to that of \cite[Algorithm 1.4.5]{St2}.
It was shown in the proof of the previous theorem that any binomial
in $I_{\cal A}$ is of the form ${\bf x}^{\mu}-{\bf x}^{\overline{\mu}}$
for some $\mu\in\cal{M_A}$.  Hence the Gr\"obner basis $G$ of $I_{\cal A}$
is of the form
$$
G=\{{\bf x}^{\mu_1}-{\bf x}^{\overline{\mu_1}},\ldots ,
{\bf x}^{\mu_r}-{\bf x}^{\overline{\mu_r}}\}.
$$

We first show that the set $\cal H$ is a generating set for $\cal{M_A}$.
Suppose not, so there exists $\mu\in\cal{M_A}$ which is not an $\mathbb{N}$-linear
combination of elements in $\cal H$.  We can choose $\mu$ so that
${\bf x}^{\mu}$ is minimal with respect to the chosen term ordering of $R$.
Since $\mu\in\cal H$, we have $\phi({\bf x}^{\mu}-{\bf x}^{\overline{\mu}})=0$,
so that ${\bf x}^{\mu}-{\bf x}^{\overline{\mu}}\in I_{\cal A}$.  Hence the
leading term of this binomial is divisible by a binomial in $G$.  That is,
there exist $i,\beta$ such that
$$
{\bf x}^\mu - {\bf x}^{\overline{\mu}} = {\bf x}^{\beta}({\bf x}^{\mu_i} - {\bf
x}^{\overline{\mu_i}}) + {\bf x}^{\overline{\mu_i}}({\bf x}^\beta - {\bf
x}^{\overline{\beta}}).
$$
Thus ${\bf x}^{\overline{\mu_i}}({\bf x}^\beta - {\bf x}^{\overline{\beta}}) \in
I_{\cal A}$. But $I_{\cal A}$ is a prime ideal, and it is immediate
from the definition of $\phi$ that it contains no
monomials, so
$ ({\bf x}^\beta - {\bf x}^{\overline{\beta}}) \in I_{\cal{A}}$.
Moreover, ${\bf x}^\beta < {\bf x}^\mu$. Hence, by the choice of $\mu$,
$\beta$ is a linear
combination of elements of ${\cal{H}}$, which implies that
$\mu = \beta + \mu_i$ is also a linear combination of elements of
${\cal{H}}$. This is a contradiction to our assumption on $\mu$.
Thus, ${\cal{H}}$ is a generating set for
${\cal{M_A}}$.

To show that $\cal H$ is minimal, suppose that
some $\mu_i$ or $\overline{\mu_i}$ is an $\mathbb{N}$-linear combination
of elements in $\cal H$.  Since the Gr\"obner basis $G$ is reduced,
the linear combination cannot contain any summands coming from
the $\mu_i,\overline{\mu_i}$.  But observe that all the vectors
${\bf e}_i+{\bf e}_j$ are symmetric, so that
$$
{\bf x}^{{\bf e}_i+{\bf e}_j}-{\bf x}^{\overline{{\bf e}_i+{\bf e}_j}}=0.
$$
A similar argument disposes of the other cases.  Thus, we have
shown that $\cal H$ is the Hilbert basis of $\cal{M_A}$.

2. $\langle {\bf x}^\mu - {\bf x}^{\overline{\mu}}\mid
\mu \in {\cal{H}} \rangle = \langle G \rangle = I_{\cal{A}}$.
\end{proof}

We can also compute the dimension of the affine variety $\vv (I_{\cal A})$
of the ideal $I_{\cal A}$.

\bt \label{dim}
The dimension of $\vv (I_{\cal A})$ is equal to $l$ if
$c_{i1}=c_{i2}$ for all $i=1,\dots,l$, and $l+1$ otherwise.
\et
\begin{proof}
According to \cite[Lemma 4.2]{St}
the dimension of $\vv (I_{\cal A})$ is equal to the number of linearly
independent column vectors in the matrix

\begin{flushleft}
$
{\cal B}=
\left(\begin{array}{ccccccccc}
1&0&\dots&0&0&\dots&0&1\\
0&1&\dots&0&0&\dots&1&0\\
\vdots&\vdots&\ddots&\vdots&\vdots&\ddots&\vdots&\vdots\\
0&0&\dots&1&1&\dots&0&0\\
 c_{11} & c_{21} & \dots  & c_{l1} &
                         c_{l2} & \dots & c_{22} & c_{12} \\
                        c_{12} & c_{22} & \dots & c_{l2} & c_{l1} & \dots &
                        c_{21} & c_{11}
\end{array} \right) \sim
$
\end{flushleft}

\begin{flushright}
$
\left(\begin{array}{ccccccccc}
1&0&\dots&0&0&\dots&0&1\\
0&1&\dots&0&0&\dots&1&0\\
\vdots&\vdots&\ddots&\vdots&\vdots&\ddots&\vdots&\vdots\\
0&0&\dots&1&1&\dots&0&0\\
 0 & 0 & \dots  & 0 &
                         c_{l2}-c_{l1} & \dots & c_{22}-c_{21} & c_{12}-c_{11}
\\
                        0 & 0 & \dots & 0 & 0 & \dots &
                        0 & 0
\end{array} \right)
$
\end{flushright}

The theorem now follows.
\end{proof}

{\it Proof of Theorems \ref{sip} and \ref{siv}}. Theorem \ref{sip}
is a corollary of Theorems \ref{iai} and \ref{hb} and Theorem \ref{siv}
follows immediately from Theorem \ref{dim}.

\medskip
To close this section, we summarize the algorithm, and
in the next section we shall provide several examples.

\smallskip
\centerline{{\bf Symmetry Component Algorithm}}

\smallskip\noindent
{\bf Input:} Two sequences of integers
$p_1,\ldots ,p_l,p_i\geq -1;q_1,\ldots ,q_l,q_i\geq 0$.
(These are the coefficient labels for a system of the form  (\ref{gs}).)

\smallskip\noindent
{\bf Output:} A finite set of generators for
the Sibirsky ideal $I_{sym}$ of (\ref{gs}), and the Hilbert basis
${\cal H}$ of the monoid $\cal M$.
\begin{enumerate}
\item Compute a reduced Gr\"obner basis $G$ for the ideal
\[
{\cal {J}} = \langle  a_{p_iq_i} - y_it_1^{p_{i}}t_2^{q_{i}}, b_{q_ip_i} 
-y_{l-i+1}t_1^{q_{{l-i+1}}}t_2^{p_{{l-i+1}}} \mid i = 1, \ldots, l \rangle
\]
in $\mathbb{C}[a,b,y_1, \dots, y_l,t_1^{\pm},t_2^{\pm}]$, 
with respect to any elimination ordering with 
\[
\{ t_1,t_2\} > \{y_1, \dots , y_d \} > \{ a_{p_1q_1} , \dots , b_{q_1p_1} \}.
\]
\item $I_{sym} =\langle G \cap \mathbb{C}[a,b]\rangle $.
\item ${\cal H}= \{ \nu, \bar{\nu} \mid [\nu] - [\bar{\nu}]
\in G\} \cup  \{ {\bf e}_i +{\bf e}_j \mid j = d-i+1; i =1,\ldots, l\}$
is a Hilbert basis for the monoid $\cal M$, where 
${\bf e}_i=(0,\dots,0,{1},0,\dots,0)$ is the $i$th basis vector.
\end{enumerate}

\section {Invariants of the Rotation Group and the Symmetry Components}

The symmetry component of the center variety of a real system in
the complex form (\ref{rgs}) was  investigated by Sibirsky
and his coworkers  \cite{Sib1,Sib2}, using invariants
of the rotation group of this system.

He showed that a monomial
$$
 a_{p_1q_1}^{\nu_{1}}a_{p_2q_2}^{\nu_{2}}\cdots
a_{p_lq_l}^{\nu_{l}}\bar a_{q_lp_l}^{\nu_{l+1}}\cdots
\bar a_{q_2p_2}^{\nu_{2l-1}}
\bar a_{q_1p_1}^{\nu_{2l}}
$$
is invariant under the action of
the rotation group
\be \label{rrot}
x\longmapsto x e^{i\varphi}
\ee
if and only if $\nu\in \mathbb{N}^{2l}$ is a solution of the Diophantine
equation
\begin{eqnarray}\label{lsib}
(p_1 -q_1)\nu_1+(p_2-q_2)\nu_2 &+&\cdots + (p_l-q_l)\nu_l\\ 
 &+& (q_l-p_l)\nu_{l+1}+\cdots +(q_1-p_1)\nu_{2l}=0, \nonumber
\end{eqnarray}
which is obtained by subtracting the second equation of system
(\ref{lk}) from the first one.
\begin{pro}
The set of all non-negative integer solutions of  equation
(\ref{lsib}) coincides with the monoid $\cal M$ defined by equation
(\ref{lk}).
\end{pro}
\begin{proof}
Obviously, every solution of (\ref{lk}) is also a solution
of (\ref{lsib}). Conversely, let $\nu$ be a solution
of (\ref{lsib}).  Then
\be \label{kr}
L^1(\nu)=L^2(\nu)=k
\ee
and
\be \label{Lpq}
L^1({\bf e}_i)+L^2({\bf e}_i)=L^1({\bf e}_{2l-i})+L^2({\bf e}_{2l-i})=p_i+q_i\ge 1,
\ee
for $i=1,\dots,l.$
Taking into account that $L(\nu)$ is a linear operator,
we conclude that $k$ on the right-hand side of (\ref{kr})
is non-negative.
\end{proof}

Thus, to find invariants of system (\ref{rgs}) under the rotation
group (\ref{rrot}) it is sufficient to find
a generating set of the monoid $\cal M$.
It also is easily seen that the monomials $[\nu],\nu \in\cal M$ are invariants
of the system (\ref{gs}) under the action of the transformation
$$
x\longmapsto x e^{i\varphi}, \ y\longmapsto y e^{-i\varphi}.
$$

Knowledge of the invariants of system (\ref{rgs}) allows one
to determine whether the corresponding vector field has an axis of symmetry.
Denote by $\tilde E$ the subset of
$E(a)$ consisting of systems (\ref{rgs})
such that the fraction $U(u,v)/V(u,v)$ is irreducible, where $U, V$ are the
functions on the right-hand
side of system (\ref{ruv}).  The following theorems are proven in \cite{Sib1}.

\begin{teo} \label{t6}
If $a^* \in  E(a)$ and ${\rm IM}[\nu]=0$
at the point $(a^*,\bar a^*)$ for all $\nu \in \cal M$, then
 the corresponding
vector  field (\ref{ruv}) has an axis of symmetry passing through the origin.
Moreover, if $a^*\in \tilde E$ then the opposite statement holds.
\end{teo}
According to this theorem, for every point in $\vv(I_{sym})$
of the form $(a^*,\bar a^*)$
the corresponding   vector field (\ref{ruv}) has an axis of symmetry
passing through the origin.
This is the reason for calling $\vv(I_{sym})$ the ``symmetry component.''
However, we have not investigated whether for
arbitrary points in $\vv(I_{sym})$ the
corresponding systems (\ref{gs}) have any kind of symmetry.

As an immediate corollary of   Theorem \ref{t6} we get
\begin{teo}
If $a^* \in  E(a)$
and
$ {\rm IM} \ [\nu]=0$
at the point $(a^*,\bar a^*)$ for all $\nu \in M$,
then  the corresponding
vector  field (\ref{ruv}) has a center at the origin.
\end{teo}
Similarly, it follows from Proposition \ref{prop1}, that
if
\be \label{ccs}
 {\rm IM}[\nu]=0
\ee
for all $\nu \in M$ at $(a^*,b^*)\in E(a,b),$
then the corresponding
system (\ref{gs}) has a center at the origin.

As mentioned above, according to \cite{Sib1,Sib2}
the monomial $[\nu]$ is an invariant of the rotation
group (\ref{rrot}) if and only if $\nu$ is a solution of (\ref{lsib}).
Proposition  2 implies that the set of solutions of
equation (\ref{lsib}) coincides with $\cal M$. The invariant
$[\nu]$ is called {\it irreducible}
if $\nu$ cannot be written in the form $\nu=\mu+\theta$,
where $\mu,\theta$ are invariants. Therefore the set $\{[\nu_i]\}$
is the set of all irreducible invariants if and only if  $\{\nu_i\}$
is the Hilbert basis of $\cal M$.

It is proven
in \cite{Sib2}, that
the maximal degree of irreducible invariants of the group
(\ref{rrot}) is less than or equal to 
\be
\label{sbound}
 2(1+\max_{(p,q)\in S}(p+q)).
\ee
Using this and equation (\ref{lk}) or  (\ref{lsib})
it is possible to find  a generating set of the monoid $\cal M$
by sorting. In particular, in \cite{Sib1}
center conditions are obtained for the cubic system
\begin{equation}\label{RCS}
i\frac{dx}{dt}=x(1-a_{10}x-a_{01}\bar x-a_{-12}x^{-1}{\bar x}^2-
a_{20}x^2-a_{11}x\bar x-
a_{02}\bar x^2-a_{-13}x^{-1}\bar x^3)
\end{equation}
which can also be obtained from Theorem \ref{tuccs} below by
taking imaginary parts.

 Theorem  \ref{hb} provides
 another way to find a Hilbert basis
of the monoid $\cal M$ and, therefore, the center conditions 
(\ref{ccs}). We apply it to find the symmetry component of
the general cubic system
\begin{eqnarray}\label{CS}
\dot x=x(1-a_{10}x-a_{01}y-a_{-12}x^{-1}y^2-a_{20}x^2-a_{11}xy-a_{02}y^
2-a_{-13}x^{-1}y^3),\\
\dot y=-y(1-b_{2,-1}x^2y^{-1}-b_{10}x-b_{01}y-b_{3,-1}x^3y^{-1}-b_{20}x
^2-b_{11}xy-b_{02}y^2),\nonumber
\end{eqnarray}
where $x,y,a_{ij},b_{ij} \in \mathbb{C}$.
\begin{teo} \label{tuccs}
The symmetry component $V_{sym}$   of the center variety
of  cubic system (\ref{CS}) is defined by the following equations:
\begin{eqnarray} \label{ucpcs}
0&=& a_{11}-b_{11}= a_{01}b_{02}b_{2,-1}-a_{-1,2}b_{10}a_{20}=a_{01}a_{02}b_{2,-1}-a_{-1,2}b_{10}b_{20} \nonumber  \\ 
&=& a_{10}^4a_{-13}-b_{3,-1}b_{01}^4=a_{10}a_{-12}b_{20}-b_{01}b_{2,-1}a_{02} =a_{10}a_{-12}b_{10}^2-a_{01}^2b_{2,-1}b_{01} \nonumber \\ 
&=& a_{20}a_{02}-b_{20}b_{02}=a_{10}^2a_{-1,2}b_{10}-a_{01}b_{2,-1}b_{01}^2=a_{10}b_{02}b_{10}-a_{01}a_{20}b_{01} \nonumber \\
&=& a_{01}^3b_{2,-1}-a_{-12}b_{10}^3 = a_{10}a_{02}b_{10}-a_{01}b_{20}b_{01} =a_{10}^3a_{-12}-b_{2,-1}b_{01}^3 \nonumber \\
&=&a_{10}a_{-1,3}b_{2,-1}-a_{-12}b_{3,-1}b_{01}= a_{20}a_{-1,3}b_{20}-a_{02}b_{3,-1}b_{02}= a_{10}^2b_{02}-a_{20}b_{01}^2 \nonumber \\
&=& a_{02}^2b_{3,-1}-a_{-13}b_{20}^2= a_{01}a_{-12}b_{3,-1}-a_{-13}b_{2,-1}b_{10}
=a_{01}^2b_{20}-a_{02}b_{10}^2 \nonumber \\
&=& a_{20}^2a_{-13}-b_{3,-1}b_{02}^2= a_{10}a_{-13}b_{20}b_{10}-a_{01}a_{02}b_{3,-1}b_{01} \nonumber \\ 
&=& a_{10}a_{20} a_{-1,3}b_{10}-a_{01} b_{3,-1} b_{02} b_{01}= a_{10}b_{02}^2b_{2,-1}-a_{-12}a_{20}^2b_{01} \nonumber \\
&=& a_{10}^2a_{02}-b_{20}b_{01}^2= a_{10}a_{02}b_{02}b_{2,-1}-a_{-12}a_{20}b_{20}b_{01} \\
&=& a_{10}a_{02}^2b_{2,-1}-a_{-12}b_{20}^2b_{01}= a_{01}^2b_{3,-1}b_{02}-a_{20}a_{-13}b_{10}^2=a_{01}^2a_{20}-b_{02}b_{10}^2 \nonumber \\
&=& a_{01}a_{-12}b_{20}^2-a_{02}^2b_{2,-1}b_{10}=a_{10}^2a_{-13}b_{20}-a_{02}b_{3,-1}b_{01}^2=a_{-12}a_{20}b_{10}-b_{02}b_{2,-1}b_{01} \nonumber \\
&=& a_{01}a_{-12}a_{20}b_{20}-a_{02}b_{02}b_{2,-1}b_{10}=a_{01}^2a_{02}b_{3,-1}-a_{-13}b_{20}b_{10}^2 \nonumber \\
&=& a_{10}^2a_{20}a_{-13}-b_{3,-1}b_{02}b_{01}^2= a_{01}a_{-12}a_{20}^2-b_{02}^2b_{2,-1}b_{10}=a_{10}a_{-13}b_{10}^3-a_{01}^3b_{3,-1}b_{01}\nonumber \\
&=&a_{10}^2a_{-13}b_{10}^2-a_{01}^2b_{3,-1}b_{01}^2= a_{10}^3a_{-13}b_{10}-a_{01}b_{3,-1}b_{01}^3=a_{01}^4b_{3,-1}-a_{-13}b_{10}^4 \nonumber \\
&=& a_{01}a_{-13}b_{20}b_{2,-1}-a_{-12}a_{02}b_{3,-1}b_{10}=a_{01}a_{20}a_{-13}b_{2,-1}-a_{-12}b_{3,-1}b_{02}b_{10} \nonumber \\
&=& a_{10}a_{-12}b_{3,-1}b_{02}-a_{20}a_{-13}b_{2,-1}b_{01}= a_{10}a_{-12}a_{02}b_{3,-1}-a_{-13}b_{20}b_{2,-1}b_{01} \nonumber \\
&=& a_{-12}^2b_{3,-1}b_{20}-a_{02}a_{-13}b_{2,-1}^2=a_{-12}^2a_{20}b_{3,-1}-a_{-13}b_{02}b_{2,-1}^2 \nonumber \\
&=& a_{10}a_{-12}^2b_{3,-1}b_{10}-a_{01}a_{-13}b_{2,-1}^2b_{01}=a_{01}^2a_{-13}b_{2,-1}^2-a_{-12}^2b_{3,-1}b_{10}^2\nonumber \\
&=& a_{-12}^2b_{20}^3-a_{02}^3b_{2,-1}^2=a_{-12}^2a_{20} b_{20}^2-a_{02}^2 b_{02}b_{2,-1}^2=a_{-12}^2a_{20}^2 b_{20}-a_{02} b_{02}^2b_{2,-1}^2 \nonumber \\
&=& a_{10}^2a_{-12}^2b_{3,-1}-a_{-13}b_{2,-1}^2b_{01}^2=a_{-12}^2a_{20}^3-b_{02}^3b_{2,-1}^2=a_{-12}^2b_{3,-1}^2b_{02}-a_{20}a_{-13}^2b_{2,-1}^2\nonumber \\
&=& a_{-12}^4b_{3,-1}^3-a_{-13}^3b_{2,-1}^4=a_{-12}^2a_{02}b_{3,-1}^2-a_{-13}^2b_{20}b_{2,-1}^2=a_{10}a_{01}-b_{10}b_{01}\nonumber \\
&=& a_{01}a_{-13}^2b_{2,-1}^3-a_{-12}^3b_{3,-1}^2b_{10}=a_{10}a_{-12}^3b_{3,-1}^2-a_{-13}^2b_{2,-1}^3b_{01}.\nonumber
\end{eqnarray}
\end{teo}
\begin{proof}
It is enough to show that the above equations (\ref{ucpcs})
form a {\it Gr\"obner basis}
of the ideal $I_{sym}$. To compute  $I_{sym}$, one can use the
{\it Symmetry Component Algorithm},
with any computer algebra system. Here we are using the specialized
system {\it Macaulay} \cite{GrSt}
which performs Gr\"obner basis calculations substantially faster than
most general purpose symbolic calculation packages.
Figure \ref{kim} shows  the Macaulay session
used to compute  $I_{sym}$ for system (\ref{CS}).
To simplify notation, we renamed the variables
$a_{ij}$ and $b_{ij}$ as follows:
$x_1 = a_{10}, x_2 = a_{01}, x_3 = a_{-12}, x_4 = a_{20}, x_5 =a_{11}, x_6 = a_{02},
x_7 = a_{-1,3}, x_8 = b_{3,-1}, x_9 = b_{20},
x_{10} = b_{11}, x_{11} = b_{02}, x_{12} = b_{2,-1}, x_{13} = b_{10},
x_{14} = b_{01}.$ 
\end{proof}

\putfig{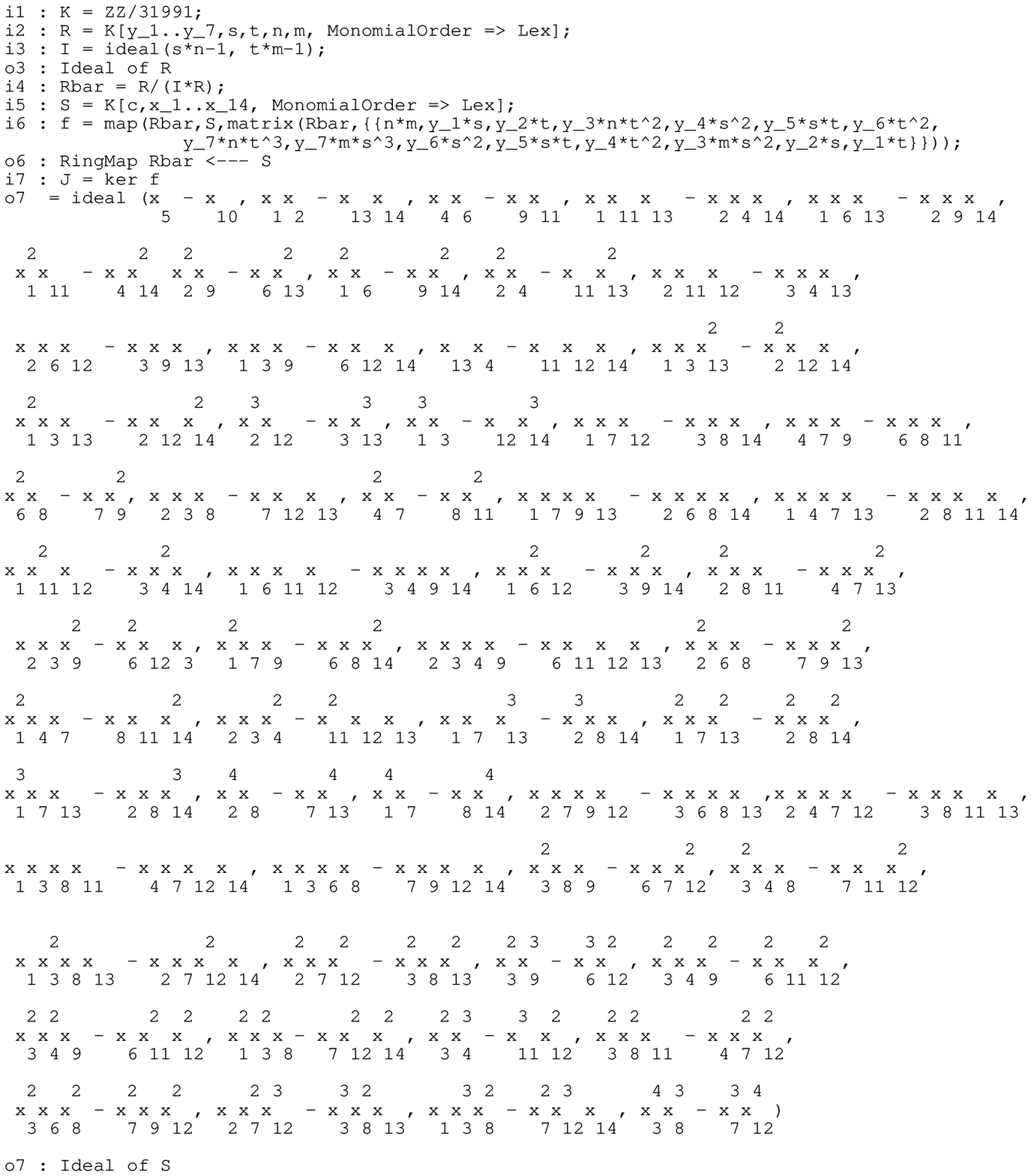}{6in}{5in}{Macaulay computation for system(\ref{CS})}{kim}

Recently some sufficient center conditions were obtained
for the real systems with homogeneous nonlinearities of
fourth and fifth degrees \cite{Chav1,Chav2}. The following are the
general complex forms of such systems.

\begin{eqnarray}\label{deg4}
\dot x=x(1-a_{30}x^3-a_{21}x^2y-a_{12}xy^2-a_{03}y^3-a_{-14}x^{-1}y^4), \\
\dot y=-y(1-b_{4,-1}x^4y^{-1}-b_{30}x^3-b_{21}x^2 y - b_{12}xy^2-b_{0,3}y^{3}),\nonumber
\end{eqnarray}
and
\begin{eqnarray}\label{deg5}
\dot x=x(1-a_{40}x^4-a_{31}x^3y-a_{22}x^{2}y^2-a_{13}xy^3-a_{04}y^4-a_{-15}x^{-1}y^5),\\
\dot y=-y(1-b_{5,-1}x^5 y^{-1}-b_{40}x^4-b_{31}x^3y-b_{22}x^2y^2-b_{13}xy^{3}-b_{04}y^4).\nonumber
\end{eqnarray}

The following theorems give some center conditions for
complex systems with homogeneous nonlinearities of fourth and
fifth degrees, respectively.

\begin{teo} \label{pf4}
The symmetry component of the center variety
of system (\ref{deg4}) is defined  by the following equations
\begin{eqnarray}
\label{ucd4}
0&=&a_{30}a_{03}- b_{30}b_{03}= a_{21}a_{12}-b_{21}b_{12}= a_{30}b_{12}^3-a_{21}^3b_{03} \nonumber \\
&=&a_{30}a_{12}b_{12}^2-a_{21}^2b_{21}b_{03}=
a_{30}a_{12}^2b_{12}-a_{21}b_{21}^2b_{03}= a_{21}a_{03}b_{21}^2-a_{12}^2-b_{30}b_{12}
\nonumber \\
&=&a_{21}^2a_{03}b_{21}-a_{12}b_{30}b_{12}^2=
a_{12}^3b_{30}-a_{03}b_{21}^3= a_{21}^3a_{03}-b_{30}b_{12}^3 \nonumber \\
&=&a_{30}a_{12}^3-b_{21}^3b_{03}=
a_{30}a_{-1,4}b_{30}b_{12}-a_{21}a_{03}b_{4,-1}b_{03} \nonumber \\
&=&a_{30}^2a_{-1,4}b_{12}-a_{21}b_{4,-1}b_{03}^2
=a_{30}a_{-1,4}b_{21}^2-a_{12}^2b_{4,-1}b_{03} \nonumber \\
&=&a_{21}a_{-1,4}b_{30}b_{21}-a_{12}a_{03}b_{4,-1}b_{12}=
a_{30}a_{21}a_{-1,4}b_{21}-a_{12}b_{4,-1}b_{12}b_{03} \nonumber \\
&=&a_{12}a_{-1,4}b_{30}^2- a_{03}^2b_{4,-1}b_{21}=
 a_{30}a_{12}a_{-1,4}b_{30}-a_{03}b_{4,-1}b_{21}b_{03} \nonumber \\
&=&a_{21}^2a_{-1,4}b_{30}-a_{03}b_{4,-1}b_{1,2}^2
=a_{21}a_{03}^2b_{4,-1}-a_{-1,4}b_{30}^2b_{12} \nonumber \\
&=&a_{30}^2a_{12}a_{-1,4}-b_{4,-1}b_{21}b_{03}^2= a_{30}a_{21}^2a_{-1,5}-
b_{4,-1}b_{12}^2b_{03} \\
&=&a_{21}a_{-1,4}b_{21}^4-a_{12}^4b_{4,-1}b_{12}=
a_{21}^2a_{-1,4}b_{21}^3-a_{12}^3b_{4,-1}b_{12}^2 \nonumber \\
&=&a_{21}^3a_{-1,4}b_{21}^2-a_{12}^2b_{4,-1}b_{12}^3
= a_{21}^4a_{-1,4}b_{21}-a_{12}b_{4,-1}b_{12}^4 \nonumber \\
&=&a_{12}^5b_{4,-1}-a_{-1,4}b_{21}^5=
a_{21}^5a_{-1,4}-b_{4,-1}b_{12}^5 \nonumber \\
&=&a_{30}a_{-1,4}^2b_{30}^2b_{21}-a_{12}a_{03}^2b_{4,-1}^2b_{03}=
a_{30}^2a_{-1,4}^2b_{30}b_{21}-a_{12}a_{03}b_{4,-1}^2b_{03}^2 \nonumber \\
&=&a_{30}^3a_{-1,4}^2b_{21}-a_{12}b_{4,-1}^2b_{03}^3
=a_{21}a_{-1,4}^2b_{30}^3-a_{03}^3b_{4,-1}^2b_{12} \nonumber \\
&=&a_{30}a_{21}a_{-1,4}^2b_{30}^2-a_{03}^2b_{4,-1}^2b_{12}b_{03}=
a_{30}^2a_{21}a_{-1,4}^2b_{30}-a_{03}b_{4,-1}^2b_{12}b_{03}^2 \nonumber \\
&=&a_{12}a_{03}^3b_{4,-1}^2-a_{-1,4}^2b_{30}^3b_{21}=
a_{30}^3a_{21}a_{-1,4}^2-b_{4,-1}^2b_{12}b_{03}^3 \nonumber \\
&=&a_{30}a_{-1,4}^3b_{30}^4-a_{03}^4b_{4,-1}^3b_{03}
=a_{30}^2a_{-1,4}^3b_{30}^3- a_{03}^3b_{4,-1}^3b_{03}^2 \nonumber \\
&=&a_{30}^3a_{-1,4}^3b_{30}^2-a_{03}^2b_{4,-1}^3b_{03}^3=
a_{30}^4a_{-1,4}^3b_{30}-a_{03}b_{4,-1}^3b_{03}^4 \nonumber \\
&=&a_{03}^5b_{4,-1}^3-a_{-1,4}^3-b_{30}^5= a_{30}^5a_{-1,4}^3-b_{4,-1}^3b_{03}^5\nonumber
\end{eqnarray}
\end{teo}
\begin{proof}
Similar to the proof of Theorem \ref{tuccs}.
To simplify notation we rename the variables
$a_{ij}$ and $b_{ij}$: $x_1 = a_{30}, x_2 = a_{21}, x_3 = a_{12}, x_4 = a_{03}, x_5 =a_{-1,4},
x_6 = b_{4,-1}, x_7 = b_{30}, x_8 = b_{21}, x_9 = b_{12}, x_{10}= b_{03}.$
Figure \ref{hh} is the {\it Macaulay} session used to find $I_{sym}$. 
\end{proof}

\putfig{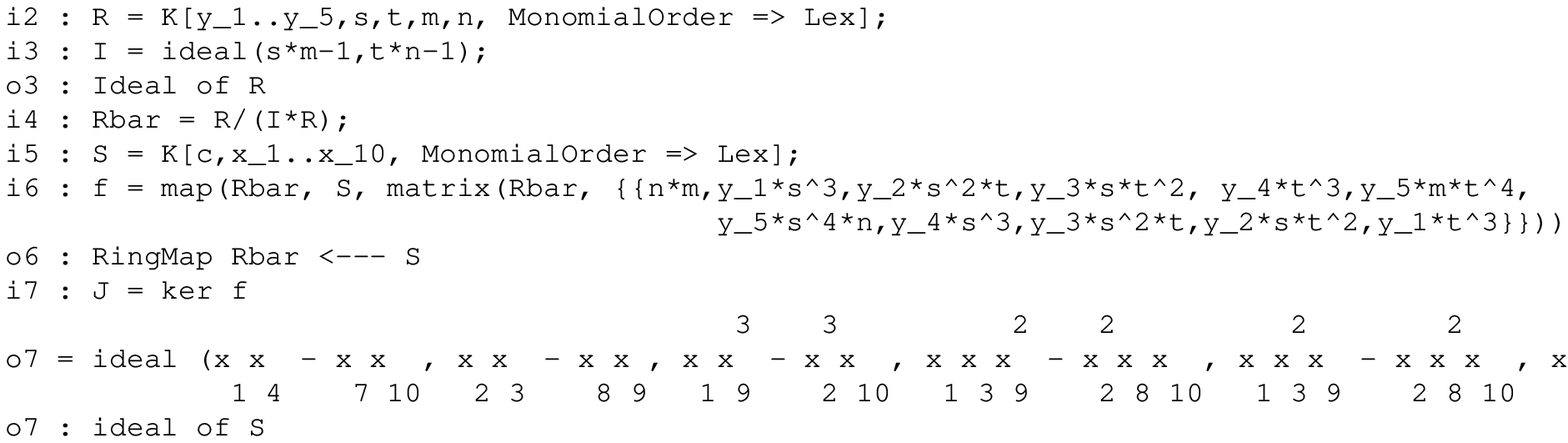}{1.6in}{4.5in}{Macaulay computation for system(\ref{deg4})}{hh}

\begin{teo} \label{fifththeorem}
The symmetry component of the center variety of 
system (\ref{deg5}) is defined  by the following equations
\begin{eqnarray}\label{ucd5}
0&=&a_{22}-b_{22}= a_{40}a_{04}-b_{40}b_{04}= a_{31}a_{13}-b_{31}b_{13}=a_{40}b_{13}^2-a_{31}^2b_{04}\nonumber \\
&=& a_{40}a_{13}b_{13}-a_{31}b_{31}b_{04}= a_{31}a_{04}b_{31}-a_{13}b_{40}b_{13}= a_{13}^2b_{40}-a_{04}b_{31}^2  \nonumber \\
&=& a_{31}^2a_{04}-b_{40}b_{13}^2= a_{40}a_{13}^2-b_{31}^2b_{04}= a_{40}a_{-1,5}b_{31}-a_{13}b_{5,-1}b_{04}\nonumber \\
&=&a_{31}a_{-1,5}b_{40}-a_{04}b_{5,-1}b_{13}= a_{13}a_{04}b_{5,-1}-a_{-1,5}b_{40}b_{31} \\
&=&a_{40}a_{31}a_{-1,5}-b_{5,-1}b_{13}b_{04}=a_{40}a_{-1,5}b_{40}a_{13}-a_{31}a_{04}b_{5,-1}b_{04} \nonumber \\
&=&a_{40}^2a_{-1,5}b_{13}-a_{31}b_{5,-1}b_{04}^2=a_{31}a_{-1,5}b_{31}^2-a_{13}^2b_{5,-1}b_{13}\nonumber \\
&=&a_{31}^2a_{-1,5}b_{31}-a_{13}b_{5,-1}b_{13}^2=a_{13}a_{-1,5}b_{40}^2-a_{04}^2b_{5,-1}b_{31} \nonumber  \\
&=&a_{40}a_{13}a_{-1,5}b_{40}-a_{04}b_{5,-1}b_{31}b_{04}=a_{31}a_{04}2b_{5,-1}-a_{-1,5}b_{40}^2b_{13}\nonumber \\
&=&a_{13}^3b_{5,-1}-a_{-1,5}b_{31}^3=a_{40}^2a_{13}a_{-1,5}-b_{5,-1}b_{31}b_{04}^2= a_{31}^3a_{-1,5}-b_{5,-1}b_{13}^3 \nonumber \\
&=& a_{40}a_{-1,5}^2b_{40}^2-a_{04}^2b_{5,-1}^2b_{04}= a_{40}^2a_{-1,5}^2b_{40}-a_{04}b_{5,-1}^2b_{04}^2 \nonumber \\
&=& a_{04}^3b_{5,-1}^2-a_{-1,5}^2b_{20}^3= a_{40}^3a_{-1,5}^2-b_{5,-1}^2b_{04}^3.\nonumber
\end{eqnarray}
\end{teo}
\begin{proof}Similar to the proof of Theorem \ref{tuccs}. \end{proof}

Thus we have presented an efficient algorithm to compute 
 the symmetry component of the center variety.
Up to now, the only known method for finding this component is due to Sibirsky.
His algorithm is as follows.
\begin{enumerate}
\item He gives the formula (\ref{sbound})  for  an upper bound for the
degrees of the irreducible invariants.
\item With this  bound one can find  all irreducible invariants by sorting.
\end{enumerate}
In \cite{Liu2} a method is given to find all ``elementary Lie invariants" 
(our Hilbert basis).  There the problem is reduced to 
finding non-negative solutions of a Diophantine equation 
similar to our equation (\ref{lsib}), which is done on
a case-by-case basis by inspection.  

Moreover neither Sibirsky nor Yi-Rong \& Ji-Bin have an analog of our Theorem 1, 
which shows that the obtained invariants generate
a prime ideal. But, as we will see in the next section, this fact is 
an important characterization of this subvariety of the center variety and 
is very helpful in the investigation 
of the cyclicity problem.

\section{Applications to Cubic Systems}

In recent years many studies have been devoted to investigating
different subfamilies of  the  cubic system (\ref{CS})
(see, e.g., \cite{Isp,DLP,DRS,Ed,RS,Liu2,ZN} and references therein).
In the case when conditions (\ref{cond}) are satisfied
 system (\ref{CS}) is equivalent to system (\ref{RCS}).
It should be mentioned that the center-focus problem is much better
investigated for the real cubic system (\ref{RCS}) than for the general
system
(\ref{CS}).

In this section we will show that   the results  obtained above, together with
additional tools from computational algebra,   give a simple  efficient  way to compute
the radical of the ideal of focus quantities and, therefore,
to solve the cyclicity problem in those cases where this ideal is radical.
As a rule it is easy to find rational parameterizations of
components of center varieties.
This fact, together with the following theorem and Theorem \ref{sip}, gives
an easy method of finding the radical of the ideal of focus quantities.
For general rational   parameterizations this is a difficult computational
problem.

\begin{teo}\label{ker}
If the  variety $\vv (J)$ of an ideal $J$
of $\mathbb{C} [x_1,\dots,x_n]$  admits a rational parameterization
\[
x_i= \frac{f_i(t_1,\dots,t_m)}{g_i(t_1,\dots,t_m)} , \;\;\; i = 1, \dots, n, \mbox{ and } 
\]
\[
\mathbb{C} [x_1,\dots,x_n]\cap \langle 1-t g, \, 
g_i(t_1,\dots,t_m)x_i - f_i(t_1,\dots,t_m) : i =1, \dots, n \rangle = J
\]
(where $g=g_1g_2\cdots g_n)$, then the ideal $J$
is a prime ideal of $\mathbb{C}[x_1,\dots,x_n]$.
\end{teo}
\begin{proof}It is sufficient to show that the ideal
\[
H =  \langle 1-t g, \, g_i(t_1,\dots,t_m)x_i - f_i(t_1,\dots,t_m) : i =1, \dots, n \rangle
\]
is prime in $\mathbb{C}[x_1,\ldots ,x_n,t_1,\ldots ,t_m,t]$. 
Consider the  ring  homomorphism
\[
\psi : \mathbb{C}[x_1,\dots,x_n,t_1,\dots,t_m,t]\longrightarrow \mathbb{C}(t_1,\dots,t_m),
\]
defined by
\[
x_i\longmapsto \frac{f_i}{g_i}, \ t_l \longmapsto t_l, \ t \longmapsto \frac{1}{g},
\]
$i=1,\dots,n, l=1,\dots,m$.
It is sufficient to prove that
$
{\rm{ker}(\psi)}=H.
$
It is clear that $H\subset  {\rm{ker}(\psi)}$.
We will show the other inclusion by induction.

Let us suppose that $h\in {\rm{ker}(\psi)} $ and $ h $ is
linear in $x$, that is,
$$
h=\sum_{i=1}^n \alpha_i(t_1,\dots,t_m,t)x_i+\alpha_0(t_1,\dots,t_m).
$$
Then
$$
\alpha_0tg=-\sum_{i=1}^n \alpha_i f_it \tilde g_i,
$$
where $\tilde g_i=g/g_i$.
Therefore
$$
h=\sum_{i=1}^n \alpha_i(x_ig_i-f_i)t \tilde g_i+(1-tg)h.
$$
 Hence, $h\in H$.

Assume now that for all polynomials of degree $k$ in $x_1,\dots,x_n$ 
and $h\in {\rm{ker}(\psi)}$, we have $ h\in H$.

Let $h\in{\rm{ker}(\psi})$ be of degree $k+1$ in $x_1,\dots,x_n$.
We can write $h$ in the form
$$
h=\sum_{i=1}^n h_i(x_i,x_{i+1},\dots,x_n,t_1,\dots,t_m,t)+
h_0(t_1,\dots,t_m,t)
$$
(here every term of $h_i$ contains $x_i$).
Consider the polynomial
$$
u=\sum_{i=1}^n \frac{h_i}{x_i} (x_ig_i-f_i)t \tilde g_i+(1-tg)h.
$$
Then $u=h+v,$ where
$$
v=-t \sum_{i=1}^n f_i\tilde g_i \frac{h_i}{x_i}-t g h_0.
$$
Since $h,u\in {\rm ker}(\psi)$ we get that  $v \in {\rm ker}(\psi)$.
Then, by the induction hypothesis, $v\in H$ and, hence, $h\in H$.
\end{proof}

We now apply the results obtained so far
to the investigation of the  cyclicity problems in some specific
cases.
First consider  the systems with homogeneous
quadratic and  cubic nonlinearities,
\be \label{cqv}
\begin{array}{l} 
\dot x=x -a_{10}x^2-a_{01}xy-a_{-12}y^2,\\
\dot y=-(y-b_{10}xy-b_{01}y^2-b_{2,-1}x^{2})
\end{array}
\ee
and 
\be \label{ch}
\begin{array}{l}
\dot x=x -a_{20} x^3-a_{11} x^2 y-a_{02}x y^2-a_{-13} y^3,\\
\dot y=-(y-b_{02} y^3-b_{11}x y^2-b_{20} x^2 y-b_{3,-1} x^3),
\end{array}
\ee
correspondingly.

In the case when 
the conditions (\ref{cond}) hold 
and the linear perturbations also are taking into account 
 the system (\ref{cqv})
corresponds to the real system on the plane $(u,v),\ x=u+iv$
\begin{equation} \label{2.1}
i\dot x=i\lambda x - x-a_{10}x^2 -a_{01}x\overline x  -
a_{-12}{\overline x}^2
\end{equation}
and the system (\ref{ch}) corresponds to
\be \label{sr}
 i\dot x=i \lambda x-  x -a_{20} x^3-a_{11} x^2 \bar x-a_{02}x \bar x^2-a_{-13}\bar x^3,
\ee
where $\lambda\in \mathbb{R}$.

For the first time the cyclicity of the origin
of the system (\ref{2.1})
was  investigated by Bautin \cite{Bau}
and of the system (\ref{sr}) by Sibirsky \cite{Sibc} 
(later on another proofs were obtained by  \.Zo\l\c{a}dek  \cite{ZQV,ZN} and Yakovenko \cite{Y}).
They proved that the following statement holds.
\begin{teo} \label{bzt}
The cyclicity of the origin of system (\ref{2.1}) 
 equals 3 \cite{Bau,ZQV,Y} and the  cyclicity of the origin of system (\ref{sr}) 
 equals 5 \cite{Sibc,ZN}.
\end{teo}

The crucial and the most difficult part  of  the proofs of this   theorem
is the following statement (see e.g. \cite{Bau,RR,Y,ZQV} for detail derivation
Theorem \ref{bzt}  from Proposition \ref{tp}).
\begin{pro} \label{tp}
1) The first three focus quantities generate the ideal of focus
   quantities of systems  (\ref{cqv}) and
 (\ref{2.1}) with $\lambda =0$.

2) The first five focus quantities generate the ideal of focus
   quantities of systems (\ref{ch}) and (\ref{sr}) with $\lambda =0$.
\end{pro}

Bautin proved Proposition \ref{tp} for the real
quadratic system (\ref{2.1})  in the Kapteyn form. The proof is
quite complicated, because the ideal of focus quantities
is not radical in this case.
A simpler proof on  Bautin's way was given by Yakovenko \cite{Y}.

\.Zo\l\c{a}dek \cite{ZQV,ZN} found a new way. He
proved the Proposition \ref{tp}
for the ring of polynomials, which are invariant under the action
of the rotation group   and  showed that the ideal of focus quantities is
radical in this ring.

We will show that with  Theorems 1 and  \ref{ker}
the proof of the  celebrated Bautin theorem, as well as the treatment 
of the cyclicity problem for the system (\ref{ch}),
becomes straightforward, using only basic knowledge of computer
algebra.

Consider first the system (\ref{ch}).
Let  $I^{(c)}=\la g_{11}, g_{22},\dots \ra $
be the ideal generated by all focus quantities of this system 
(the so-called {\it Bautin ideal})  and let $I^{(c)}_k=\la g_{11}, g_{22},\dots, g_{kk} \ra $ be the ideal generated
by the first $k$ focus quantities.

Computing   the first five  focus quantities by means of the
algorithm given in \cite{R93}, and then reducing them, we find
that
\begin{eqnarray*}
g_{11}&=&a_{11}-b_{11};\\
g_{22}&=&a_{20} a_{02} - b_{02} b_{20};\\
g_{33}&=&(3 a_{20}^2 a_{-13} + 8 a_{20} a_{-13} b_{20} + 3 a_{02}^2 b_{3,-1} \\
      & & \mbox{} - 8 a_{02} b_{02} b_{3,-1} - 3 a_{-13} b_{20}^2
          - 3 b_{02}^2 b_{3,-1})/8;\\
g_{44}&=&( - 9 a_{20}^2 a_{-13} b_{11} + a_{11} a_{-13} b_{20}^2 + 9 a_{11}
           b_{02}^2 b_{3,-1} - a_{02}^2 b_{11} b_{3,-1})/16;\\
g_{55}&=&( - 9 a_{20}^2 a_{-13} b_{02} b_{20} + a_{20} a_{02} a_{-13} b_{20}^2
            + 9 a_{20} a_{02} b_{02}^2 b_{3,-1} +\\
      & & 18 a_{20} a_{-13}^2 b_{20} b_{3,-1} + 6 a_{02}^2 a_{-13} b_{3,-1}^2
	- a_{02}^2 + \\
      & & b_{02} b_{20} b_{3,-1} - 18 a_{02} a_{-13}
           b_{02} b_{3,-1}^2 - 6 a_{-13}^2 b_{20}^2 b_{3,-1})/36.
\end{eqnarray*}
Furthermore, the variety of the ideal $I$ coincides with the variety of
the first five focus quantities. More precisely, the following result holds.
\bt \label{ccch}
The center variety $\vv (I^{(c)}) $ of the system  (\ref{ch}) consists of
the three irreducible components:
\be \label{dec}
\vv (I^{(c)})=\vv (I^{(c)}_5)=\vv (J_1)\cup \vv (J_2) \cup \vv (J_3),
\ee
where
 $ J_1=\la   a_{11}-b_{11}, 3 a_{20}-b_{20}, 3 b_{02} -a_{02}\ra ,\
  J_2=\la a_{11}, b_{11}, a_{20}+3 b_{20}, b_{02}+3 a_{02},  a_{-13}
b_{3,-1}-4 a_{02} b_{20}\ra $ and
  $ J_3=\la a_{20}^2 a_{-13} - b_{3,-1} b_{02}^2,
a_{20} a_{02} - b_{20} b_{02}, $ $
a_{20} a_{-13} b_{20} - a_{02} b_{3,-1} b_{02},
a_{11} - b_{11},
a_{02}^2 b_{3,-1} - a_{-13} b_{20}^2 \ra .
$
\et
\begin{proof}
It is easy to check (using, for example, the radical membership test, see e.g.
\cite{Cox})
that
$
\vv (I^{(c)}_5)=\vv (J_1)\cup \vv (J_2) \cup \vv (J_3).
$
Thus we only have to show that $g_{mm}|_{\vv (J_k)}\equiv 0$ for all $m>0$ and
$k=1,2,3.$

Indeed, the systems corresponding to the points of $\vv (J_1)$
are Hamiltonian.

For the variety  $\vv (J_2)$
in the case $
a_{-13} = (4 a_{02} b_{20})/b_{3,-1}$
one can easily find an invariant conic $l_1$ and invariant cubic $l_2$
and check that  the system has the first integral
\be \label{min}
\Phi=l_1^3 l_2^{-2},
\ee
defined on $\vv (J_2)\setminus \vv (H)$, where
$H=\la  b_{20} b_{3,-1}^2,  b_{3,-1}^2 (a_{02} b_{3,-1}^2 + 4 b_{20}^3)\ra.$
Computing the quotient of the ideals (e.g. by means of  the algorithm from
\cite{Cox})
we get $J_2:H=J_2$.
Therefore,
$$
  \overline{\vv (J_2)\setminus \vv (H)}=\vv (J_2).
$$
This implies that the system (\ref{ch}) has a center on the whole component
$\vv (J_2).$

Finally, according to Theorem \ref{tuccs},  $\vv (J_3)$ is the symmetry component of the center variety.

The irreducibility of the  components $\vv (J_1)$ and $\vv (J_2)$
is obvious, and $\vv (J_3)$ is irreducible because due to  Theorems \ref{sip}
and \ref{tuccs}  the ideal $J_3$ is prime.
\end{proof}

A very similar theorem is proven also in \cite{Liu2} (and in  \cite{Isp1,ZN} for the 
real case).  However, our theorem 
also shows that (\ref{dec}) is the irreducible decomposition of 
the center variety.

It is worthwhile  to mention that, if one presents 
the solution of the center problem for a system of the type (\ref{gs}) 
writing out the irreducible decomposition 
of a center variety, then 
the answer is  unique. 
Otherwise it can happen that the center conditions for the 
same system obtained by different authors can look very different.
In fact   
the center conditions for  system (\ref{ch})
were first obtained by Sadovsky \cite{SadC}. He found eight such
sets but they look very different from the ones presented in Theorem 
\ref{ccch}.  However computing the Zariski closure of  the sets given by Sadovsky
and taking intersections of the corresponding ideals using
standard algorithms  from computational algebra (see below)
we get the ideal $I^{(c)}_5.$
This implies that his center conditions coincide with the three components
given in Theorem \ref{ccch}.

Recall that we can use the following algorithm to compute
the intersection of two ideals \cite{Cox}. Let $I=\la f_1,\dots, f_r\ra$
and $J=\la h_1,\dots,h_s\ra$ be ideals in $k[x_1,\dots,x_n]$.
Compute a Gr\"obner basis for the ideal
$$
\la tf_1,\dots,tf_r, (1-t)h_1,\dots,(1-t) h_s\ra \subset
k[t,x_1,\dots,x_n]
$$
using a lexicographic term order
with $t$ greater than the $x_i$. Those elements of this  basis  which  do
not contain the
variable $t$ will form a basis of $I\cap J.$

For the case of quadratic system 
one can easily compute the three first focus quantities and get
 
\vspace{.2cm}
\noindent
$
g^{(q)}_{11}=  a_{10} a_{01}-b_{10}b_{01}, \vspace{.1cm} \newline
g^{(q)}_{22}= a_{10} a_{-12}{b_{10}}^2 - a_{01}^2b_{01} {b_{2,-1}} - 
  \frac{2}3(a_{-12} b_{10}^3 - a_{01}^3 b_{2,-1})  - \frac{2}3(a_{01} b_{01}^2 b_{2,-1}- a_{10}^2 a_{-12} b_{10}),\vspace{.1cm} \newline 
g^{(q)}_{33}= \frac 5{12} (a_{01}^2 a_{10} b_{01}^2 b_{2,-1}-a_{10}^2 a_{-12} b_{01} b_{10}^2) - \frac 5{48} (a_{01}^3 b_{01} b_{10} b_{2,-1}-a_{01} a_{10} a_{-12} b_{10}^3)\vspace{.1cm} \newline\mbox{ }\hspace{.5cm} + \frac 58 (a_{01} a_{-12} b_{01}^2 b_{2,-1}^2-a_{10}^2 a_{-12}^2 b_{10} b_{2,-1}) - \frac  5{16}(a_{01}^2 a_{-12} b_{01} b_{2,-1}^2-a_{10} a_{-12}^2 b_{10}^2 b_{2,-1}). \vspace{.2cm}
$

\noindent
Here and below the upper index $(q)$ means that we are speaking about 
the quadratic system.

Using the equations above and Theorems 1 and \ref{ker} it is easy to see
(similarly to the case of the system (\ref{ch})) that the following theorem holds
(see \cite{Dul,Isp1,Sib1,ZQV} for details).
\begin{teo}\label{tz} 
The center variety  of the system (\ref{cqv}) consists of four irreducible
components:
\begin{enumerate}
\item $\vv (J_1^{(q)})$, where $J_1^{(q)}=\la 2 a_{10} -b_{10}, 2 b_{01}- a_{01}\ra,$
\item $\vv (J_2^{(q)}),$ where $J_2^{(q)}=\la a_{01},b_{10}\ra,$
\item $\vv (J_3^{(q)})$, where $J_3^{(q)}=\la 2 a_{01}+b_{01}, a_{10}+2 b_{10},
a_{01} b_{10} -a_{-12} b_{2,-1}\ra,$
\item $\vv (J_4^{(q)})=\la f_1, f_2, f_3, f_4, f_5\ra,$
where
$f_1=a_{01}^3 b_{2,-1}-a_{-12} b_{10}^3$,
$f_2=a_{10} a_{01} - b_{01} b_{10}$,
$f_3=a_{10}^3 a_{-12}- b_{2,-1} b_{01}^3$,
$f_4=a_{10} a_{-12} b_{10}^2- a_{01}^2 b_{2,-1} b_{01}$,
$f_5=a_{10}^2 a_{-12} b_{10} - a_{01} b_{2,-1} b_{01}^2.$
\end{enumerate}
\end{teo}

To our knowledge the notion of the center variety 
 was introduced in the literature  
on the center problem very  recently by \.Zo\l\c{a}dek in \cite{Z},
where he also gives Theorem \ref{tz}. But  he   writes out 
the fourth component in the form 
$$
f_1=f_2=0.
$$
It is not precise because, using Theorem \ref{ker}, one can see that
$$
\vv (\la f_1, f_2 \ra)= \vv (J_4^{(q)})  \cup \vv (J_2^{(q)}).
$$

\begin{lem}
The ideals $I^{(q)}$ and   $I^{(c)}$   are radical ideals in 
$\mathbb{C}[a_{10},a_{01},a_{-12},b_{2,-1},b_{10},b_{01}]$ and 
$\mathbb{C}[a_{20},a_{11},\dots,b_{02}]$, respectively.
\end{lem}
\begin{proof}
Computing the intersection of the ideals $J_k$ we find
$$
I^{(c)}_5=J_1\cap J_2 \cap J_3.
$$
Hence  $I^{(c)}_5$ is radical because due to 
Theorem {\ref{ker}}  $J_1, J_2$ are
prime, and according to Theorems \ref{sip} and \ref{tuccs}  $J_3$ is 
  prime as well.

Similarly, for quadratic system  we easily check  that, 
for  $1 \leq i \leq 4$, the ideals $J_i^{(q)}$ are prime and
\be \label{inters}
I^{(q)}_3=\cap^4_{i=1} J_i^{(q)}.
\ee
This yields that the ideal of focus quantities of quadratic system  is
a radical ideal.
\end{proof}

\begin{proof}(Proposition \ref{tp})
$\vv (I^{(q)})=\vv (I_3^{(q)}) $ 
($\vv (I^{(c)})=\vv (I_5^{(c)})$) 
 and  the ideal $I_3^{(q)}$ ($I_5^{(c)}$)
is a radical ideal. Therefore $I^{(q)}=I_3^{(q)}$ ($I^{(c)}=I_5^{(c)}).$
\end{proof}

In \cite{JLR}, we study the following system.
\be \label{some}
\begin{array}{l}
\dot x=x -a_{10} x^2-a_{01} x y-a_{-13} y^3,\\
\dot y=-(y-b_{01} y^2-b_{10} x y-b_{3,-1} x^3).
\end{array}
\ee
Using the methods of this paper, we compute the ideal of
the center variety of (\ref{some}) which  turns not to be radical.
Moreover, we prove the following theorem.
\begin{theorem}
The cyclicity of the origin of the system
$$
i{\dot x} = x- a_{10}x^2 - a_{01}x \bar x- a_{-13}\bar{x}^{3},
$$
is less than or equal to 5.
\end{theorem}

\section{Acknowledgments}
This work was begun during the 1999 Summer Conference
of the Rocky Mountain Mathematics Consortium on
Computational Algebra and its Applications at the University of
Wyoming.
The authors thank the Consortium for support, and the University of
Wyoming for its hospitality.
The first and second authors thank the Physical Science Laboratory
at New Mexico State University for partial support of this
work.  The third author acknowledges  support of the research through a
grant of the Ministry of
Science
and Technology of the Republic of Slovenia and  the Abdus Salam   ICTP
(Trieste) Joint Programme, the grant of  the
Foundation of Fundamental Research of the Republic of Belarus
and the sponsorship of {\it Nova Kreditna Banka Maribor}.
The authors also thank Karin Gatermann, Michael Singer, and Bernd
Sturmfels for helpful comments.


\end{document}